\documentclass[11pt,reqno]{amsart}
\usepackage{amsmath,amssymb}

 \makeatletter
 \evensidemargin  \oddsidemargin
 \makeatother

\begin{document}
\thispagestyle{empty}
 \setcounter{page}{1}

\begin{center}
{\large\bf Stability of generalized  mixed type
additive-quadratic-cubic functional equation in non-Archimedean
spaces

\vskip.15in

{\bf M. Eshaghi Gordji} \\[2mm]

{\footnotesize Department of Mathematics,
Semnan University,\\ P. O. Box 35195-363, Semnan, Iran\\
[-1mm] e-mail: {\tt madjid.eshaghi@gmail.com }}

{\bf  M. Bavand Savadkouhi} \\[2mm]

{\footnotesize Department of Mathematics,
Semnan University,\\ P. O. Box 35195-363, Semnan, Iran\\
[-1mm] e-mail: {\tt  bavand.m@gmail.com }}

{\bf \bf Th. M. Rassias} \\[2mm]
{\footnotesize Department of Mathematics, National Technical
University
\\of Athens,
 Zografou, Campus 15780 Athens, Greece\\
[-1mm]e-mail: {\tt trassias@math.ntua.gr}}

}
\end{center}
\vskip 5mm

 \noindent{\footnotesize{\bf Abstract.}
In this paper, we prove generalized Hyres--Ulam--Rassias stability
of the mixed type additive, quadratic and cubic functional
equation
$$f(x+ky)+f(x-ky)=k^2f(x+y)+k^2f(x-y)+2(1-k^2)f(x)$$ for fixed
integers $k$ with $k\neq0,\pm1$ in non-Archimedean spaces.
 \vskip.10in

  \newtheorem{df}{Definition}[section]
  \newtheorem{rk}[df]{Remark}
   \newtheorem{lem}[df]{Lemma}
   \newtheorem{thm}[df]{Theorem}
   \newtheorem{pro}[df]{Proposition}
   \newtheorem{cor}[df]{Corollary}
   \newtheorem{ex}[df]{Example}

 \setcounter{section}{0}
 \numberwithin{equation}{section}

\vskip .2in

\begin{center}
\section{Introduction}
\end{center}

We say that a functional equation $(\xi)$ is stable if any function
$g$ satisfying the equation $(\xi)$ {\it approximately} is near to
true solution of $(\xi).$  We say that a functional equation $(\xi)$
is
 {\it superstable} if every approximately solution is an exact solution
 of the equation $(\xi)$  \cite{J,R4,R2}.\\
The stability problem of functional equations originated from a
question of Ulam \cite{Ul} in 1940, concerning the stability of
group homomorphisms. Let $(G_1,.)$ be a group and let $(G_2,*)$ be a
metric group with the metric $d(.,.).$ Given $\epsilon >0$, does
there exist a $\delta> 0$, such that if a mapping
$h:G_1\longrightarrow G_2$ satisfies the inequality
$d(h(x.y),h(x)*h(y)) <\delta$ for all $x,y\in G_1$, then there
exists a homomorphism $H:G_1\longrightarrow G_2$ with
$d(h(x),H(x))<\epsilon$ for all $x\in G_1?$ In the other words,
under what condition does there exist a homomorphism near an
approximate homomorphism? The concept of stability for functional
equation arises when we replace the functional equation by an
inequality which acts as a perturbation of the equation. In 1941, D.
H. Hyers \cite{Hy} gave the first affirmative  answer to the
question of Ulam for Banach spaces. Let $f:{E}\longrightarrow{E'}$
be a mapping between Banach spaces such that
$$\|f(x+y)-f(x)-f(y)\|\leq \delta $$
for all $x,y\in E,$ and for some $\delta>0.$ Then there exists a
unique additive mapping $T:{E}\longrightarrow{E'}$ such that
$$\|f(x)-T(x)\|\leq \delta$$
for all $x\in E.$ Moreover if $f(tx)$ is continuous in
$t\in\mathbb{R}$  for each fixed $x\in E,$ then $T$ is linear.
 In $1978,$ Th. M. Rassias \cite{Ra} provided a
generalization of Hyers' Theorem which allows the Cauchy difference
to be unbounded.\\ The functional equation
$$f(x+y)+f(x-y)=2f(x)+2f(y),\eqno(1.1)$$
is related to symmetric bi-additive function
\cite{Ac-Dh},\cite{Am},\cite{Jo-Ne} and \cite{Ka}. It is natural
that this equation is called a quadratic functional equation. In
particular, every solution of the quadratic equation $(1.1)$ is
said to be a quadratic function. It is well known that a function
$f$ between real vector spaces is quadratic if and only if there
exits a unique symmetric bi-additive function $B$ such that
$f(x)=B(x,x)$ for all $x$ (see \cite{Ac-Dh},\cite{Ka}). The
bi-additive function $B$ is given by
$$B(x,y)=\frac{1}{4}(f(x+y)-f(x-y)).$$
A Hyers--Ulam--Rassias stability problem for the quadratic
functional equation $(1.1)$ was proved by Skof for functions
$f:A\longrightarrow B$, where $A$ is normed space and $B$ Banach
space (see \cite{Sk}). Cholewa \cite{Ch} noticed that the Theorem of
Skof is still true if relevant domain $A$ is replaced an abelian
group (see also  \cite{Cz} and  \cite{Gr}).\\ Jun and Kim
\cite{Ju1-Ki1} introduced the following cubic functional equation
$$f(2x+y)+f(2x-y)=2f(x+y)+2f(x-y)+12f(x), \eqno(1.2)$$ and they
established the general solution and the generalized
Hyers--Ulam--Rassias stability for the  functional equation
$(1.2).$ The $f(x)=x^3$ satisfies the functional equation $(1.2),$
which is called a cubic functional equation. Every solution of the
cubic functional equation is said to be a cubic function. Jun and
Kim proved that  a function $f$ between real vector spaces X and Y
is a solution of $(1.2)$ if and only if there exits a unique
function $C:X\times X\times X\longrightarrow Y$ such that
$f(x)=C(x,x,x)$ for all $x\in X,$ and $C$ is
symmetric for each fixed one variable and is additive for fixed two variables.\\

The stability problems of several functional equations have been
extensively investigated by a number of authors and there are many
interesting results concerning this problem (see
\cite{Ga,G,Gr,H-I-R,I-R,R1,R2,R4} and \cite{T}).

By a non-Archimedean field we mean a field $K$ equipped with a
function (valuation) $|.|$ from $K$ into $[0,\infty)$ such that
$|r|= 0$ if and only if $r=0,$ $|rs|=|r||s|,$ and $|r+s| \leq \max
\{|r|, |s|\}$ for all $r,s \in K.$ Clearly $|1|=|-1|=1$ and $|n|
\leq 1$ for all $n \in N.$\\
\begin{df}Let $X$ be a vector space over a scalar field $K$ with a
non-Archimedean non-trivial valuation $|.|.$ A function $\|.\|:X
\to \Bbb R$ is a non-Archimedean norm (valuation) if it satisfies
the following conditions:\\ $(i)$ $\|x\|=0$ if and only if
$x=0;$\\
$(ii)$ $\|rx\|=|r| \|x\| $ for all $r \in  K,$ $x \in X;$\\
$(iii)$ the strong triangle inequality (ultrametric); namely,
$$\|x+y\| \leq \max \{ \|x\|,\|y\| \}\hspace{1 cm} (x,y \in X)$$
Then $(X,\|.\|)$ is called a non-Archimedean space.
\end{df}
Due to the fact that $$\|x_n-x_m\| \leq \max \{\|x_{j+1}-x_{j}\| :
m \leq j \leq n-1 \} \hspace{1cm} (n> m)$$ a sequence $\{x_n\}$ is
Cauchy if and only if $\{x_{n+1}-x_n \}$ converges to zero in a
non-Archimedean space. By a complete non-Archimedean space we mean
one in which every Cauchy sequence is convergent.\\
M. S. Moslehian and Th. M. Rassias \cite{Mo-Ra} proved the
generalized Hyers--Ulam stability of the Cauchy functional equation
and the quadratic functional equation in non-Archimedean spaces. M.
Eshaghi Gordji and M. Bavand Savadkouhi \cite{Es1-Ba1}, have
obtained the generalized Hyers--Ulam--Rassias stability for cubic
and quartic functional equation in non-Archimedean spaces.\\

Recently, M. Eshaghi Gordji and H. Khodaei \cite{E-K}, investigated
the solution and stability of the generalized mixed type cubic,
quadratic and additive functional
equation$$f(x+ky)+f(x-ky)=k^2f(x+y)+k^2f(x-y)+2(1-k^2)f(x)\eqno(1.3)$$
for fixed integers $k$ with $k\neq0,\pm1$  in quasi--Banach spaces.
We only mention here the papers \cite{E-Z-B-R,G1} and \cite{G2}
concerning the stability of the mixed type functional equations. In
this paper, we prove  the stability of functional equation (1.3)
 in non-Archimedean
space.

\vskip.2in
\section{Stability}
Throughout this section, we assume that $G$ is an additive group and
$X$ is a complete non-Archimedean space. Given $f:G \to X,$ we
define the difference operator
$$Df(x,y)=f(x+ky)+f(x-ky)-k^2f(x+y)-k^2f(x-y)-2(1-k^2)f(x)$$
for fixed integers $k$ with $k\neq0,\pm1$ and for all $x,y \in G.$
We consider the following function inequality: $$\|Df(x,y)\| \leq
\varphi(x,y)$$ for an upper bound: $\varphi:G \times G \to
[0,\infty).$

\begin{thm}
Let $\varphi:G \times G \to [0,\infty)$ be a function such that
$$\lim_{n \to \infty} \frac{\varphi(k^nx,k^ny)}{|k|^{2n}}=0\eqno(2.1)$$
$$\lim_{n \to \infty} \frac{1}{|2.k^{2n}|}\varphi(0,k^{n-1}x)=0 \eqno(2.2)$$ for all $x,y \in G$ and let for each $x \in G$ the
limit
$$\lim_{n \to \infty} \max \{\frac{1}{|k^{2j}|}\varphi(0,k^jx) :~ 0 \leq j < n \},\eqno(2.3)$$ denoted by $\tilde{\varphi}_Q(x),$
exist. Suppose that $f:G \to X$ is an even function satisfying
$$\|Df(x,y)\| \leq \varphi(x,y) \eqno(2.4)$$
for all $x,y \in G.$ Then there exist a quadratic function $Q:G
\to X$ such that
$$\|Q(x)-f(x)\| \leq \frac{1}{|2.k^2|} \tilde{\varphi}_Q(x) \eqno(2.5)$$
for all $x \in G.$ Moreover, if
\begin{align*}
\lim_{i \to \infty} \lim_{n \to \infty}
 \max \{\max \{\frac{1}{|k^{2j}|}\varphi(0,k^jx)
\}:~i \leq j < n+i \}=0,
\end{align*}
then $Q$ is the unique quadratic function satisfying $(2.5).$
\end{thm}
\begin{proof}
By putting $x=0$ in $(2.4),$ we get
$$\|2f(ky)-2k^2f(y)\| \leq \varphi(0,y) \eqno(2.6)$$ for all $y\in G.$
If we replace $y$ in $(2.6)$ by $x,$ and divide both sides of
$(2.6)$ by $2k^2,$ we get
$$\|\frac{f(kx)}{k^2}-f(x)\| \leq \frac{1}{|2k^2|} \varphi(0,x) \eqno(2.7)$$ for all $x\in G.$
Replacing $x$ by $k^{n-1}x$ in $(2.7),$ we get
$$\|\frac{f(k^nx)}{k^{2n}}-\frac{f(k^{n-1}x)}{k^{2(n-1)}}\| \leq
\frac{1}{|2.k^{2n}|}\varphi(0,k^{n-1}x) \eqno(2.8)$$ for all $x
\in G.$ It follows from $(2.2)$ and $(2.8)$ that the sequence
$\{\frac{f(k^nx)}{k^{2n}}\}$ is Cauchy. Since $X$ is complete, we
conclude that $\{\frac{f(k^nx)}{k^{2n}}\}$ is convergent. Set
$Q(x):=\lim_{n \to \infty} \frac{f(k^nx)}{k^{2n}}.$\\
Using induction one can show that
$$\|\frac{f(k^nx)}{k^{2n}}-f(x)\| \leq  \frac{1}{|2.k^2|} \max \{\frac{1}{|k^{2i}|}\varphi(0,k^ix) :~ 0 \leq i < n \} \eqno(2.9)$$ for all $n \in \Bbb N$ and
all $x \in G.$ By taking $n$ to approach infinity in $(2.9)$ and
using $(2.3)$ one obtains $(2.5).$ By $(2.1)$ and $(2.4),$ we get
$$\|DQ(x,y)\|=\lim_{n \to \infty}\frac{1}{|k^{2n}|}\|f(k^nx,k^ny)\| \leq \lim_{n \to \infty} \frac{\varphi(k^nx,k^ny)}{|k|^{2n}}=0$$
for all $x,y \in G.$ Therefore the function $Q:G \to X$ satisfies
$(1.3).$ If $Q^{'}$ is another quadratic function satisfying
$(2.5),$ then
\begin{align*}
\|Q(x)-Q^{'}(x)\|&=\lim_{i \to
\infty}|k|^{-2i}\|Q(k^ix)-Q^{'}(k^ix)\|\\
& \leq \lim_{i \to \infty} |k|^{-2i} \max \{~
\|Q(k^ix)-f(k^ix)\|,\|f(k^ix)-Q^{'}(k^ix)\|~ \}\\
&\leq \frac{1}{|2.k^2|} \lim_{i \to \infty} \lim_{n \to \infty}
 \max \{\max \{\frac{1}{|k^{2j}|}\varphi(0,k^jx)
\}:~i \leq j < n+i \}\\
&=0.
\end{align*}
for all $x \in G.$ Therefore $Q=Q^{'}.$ This completes the proof
of the uniqueness of $Q.$
\end{proof}

\begin{thm}
Let $\varphi:G \times G \to [0,\infty)$ be a function such that
$$\lim_{n \to \infty} \frac{1}{|2^n|}\max\{\varphi(2^{n+1}x,2^{n+1}y),|8|\varphi(2^nx,2^ny)\}=0\eqno(2.10)$$
\begin{align*}
\lim_{n \to \infty} \frac{1}{|2^n.k^2(k^2-1)|}
&\max\{~\max\{\max\{|2(k^2-1)|\varphi(2^{n-1}x,2^{n-1}x),|k^2|\varphi(2^nx,2^{n-1}x)\}\\&,\max\{\varphi(2^{n-1}x,2^nx)
,\max\{\varphi(2^{n-1}(k+1)x,2^{n-1}x),\varphi(2^{n-1}(k-1)x,2^{n-1}x)\}~\}\}\\
&,\max\{\max\{\varphi(2^{n-1}x,2^{n-1}x),|k^2|\varphi(2^nx,2^nx)\}\\&,\max\{\max\{|2(k^2-1)|\varphi(2^{n-1}x,2^nx)
,\varphi(2^{n-1}x,3.2^{n-1}x)\}\\
&,\max\{\varphi(2^{n-1}(2k+1)x,2^{n-1}x),\varphi(2^{n-1}(2k-1)x,2^{n-1}x)\}\}~\}\}=0\\
&\hspace{9cm}(2.11)
\end{align*}
for all $x,y \in G$ and let for each $x \in G$ the limit
\begin{align*}
\max
&\{\frac{1}{|2^i|}\max\{~\max\{\max\{|2(k^2-1)|\varphi(2^{i-1}x,2^{i-1}x),|k^2|\varphi(2^ix,2^{i-1}x)\}\\&,\max\{\varphi(2^{i-1}x,2^ix)
,\max\{\varphi(2^{i-1}(k+1)x,2^{i-1}x),\varphi(2^{i-1}(k-1)x,2^{i-1}x)\}~\}\}\\
&,\max\{\max\{\varphi(2^{i-1}x,2^{i-1}x),|k^2|\varphi(2^ix,2^ix)\}\\&,\max\{\max\{|2(k^2-1)|\varphi(2^{i-1}x,2^ix)
,\varphi(2^{i-1}x,3.2^{i-1}x)\}\\
&,\max\{\varphi(2^{i-1}(2k+1)x,2^{i-1}x),\varphi(2^{i-1}(2k-1)x,2^{i-1}x)\}\}~\}\}:~
0 \leq i < n \},\hspace{1cm}(2.12)
\end{align*}
denoted by $\tilde{\varphi}_A(x),$ exist. Suppose that $f:G \to X$
is an odd function satisfying
$$\|Df(x,y)\| \leq \varphi(x,y) \eqno(2.13)$$
for all $x,y \in G.$ Then there exist an additive function $A:G
\to X$ such that
$$\|f(2x)-8f(x)-A(x)\| \leq \frac{1}{|2.k^2(k^2-1)|} \tilde{\varphi}_A(x) \eqno(2.14)$$
for all $x \in G.$ Moreover,  if
\begin{align*}
\lim_{i \to \infty} \lim_{n \to \infty} \max
&\{\frac{1}{|2^j|}\max\{~\max\{\max\{|2(k^2-1)|\varphi(2^{j-1}x,2^{j-1}x),|k^2|\varphi(2^jx,2^{j-1}x)\}\\&,\max\{\varphi(2^{j-1}x,2^jx)
,\max\{\varphi(2^{j-1}(k+1)x,2^{j-1}x),\varphi(2^{j-1}(k-1)x,2^{j-1}x)\}~\}\}\\
&,\max\{\max\{\varphi(2^{j-1}x,2^{j-1}x),|k^2|\varphi(2^jx,2^jx)\}\\&,\max\{\max\{|2(k^2-1)|\varphi(2^{j-1}x,2^jx)
,\varphi(2^{j-1}x,3.2^{j-1}x)\}\\
&,\max\{\varphi(2^{j-1}(2k+1)x,2^{j-1}x),\varphi(2^{j-1}(2k-1)x,2^{j-1}x)\}\}~\}\}:~
i \leq j < n+i \}\\
&=0,
\end{align*}
then $A$ is the unique additive function satisfying $(2.14).$
\end{thm}
\begin{proof}
It follows from $(2.13)$ and using oddness of $f$ that
$$\|f(ky+x)-f(ky-x)-k^2f(x+y)-k^2f(x-y)+2(k^2-1)f(x)\| \leq \varphi(x,y) \eqno(2.15)$$
for all $x,y \in G.$ Putting $y=x$ in $(2.15),$ we have
$$\|f((k+1)x)-f((k-1)x)-k^2f(2x)+2(k^2-1)f(x)\|\leq \varphi(x,x) \eqno(2.16)$$
for all $x\in G.$ It follows from $(2.16)$ that
$$\|f(2(k+1)x)-f(2(k-1)x)-k^2f(4x)+2(k^2-1)f(2x)\|\leq
\varphi(2x,2x) \eqno(2.17)$$ for all $x\in G.$ Replacing $x$ and
$y$ by $2x$ and $x$ in $(2.15),$ respectively, we get
$$\|f((k+2)x)-f((k-2)x)-k^2f(3x)-k^2f(x)+2(k^2-1)f(2x)\|\leq
\varphi(2x,x)  \eqno(2.18)$$ for all $x\in G.$ Setting $y=2x$ in
$(2.15),$ gives
$$\|f((2k+1)x)-f((2k-1)x)-k^2f(3x)-k^2f(-x)+2(k^2-1)f(x)\|\leq
\varphi(x,2x)  \eqno(2.19)$$
 for all $x\in G.$
Putting $y=3x$ in $(2.15),$ we obtain
$$\|f((3k+1)x)-f((3k-1)x)-k^2f(4x)-k^2f(-2x)+2(k^2-1)f(x)\|\leq
\varphi(x,3x)  \eqno(2.20)$$ for all $x\in G.$ Replacing $x$ and
$y$ by $(k+1)x$ and $x$ in $(2.15),$ respectively, we get
\begin{align*}
\|f((2k+1)x)-f(-x)-k^2f((k+2)x)-k^2f(kx)+2(k^2-1)&f((k+1)x)\|\\&\leq
\varphi((k+1)x,x) \hspace {.7cm}(2.21)
\end{align*}
for all $x\in G.$ Replacing $x$ and $y$ by $(k-1)x$ and $x$ in
$(2.15),$ respectively, one gets
\begin{align*}
\|f((2k-1)x)-f(x)-k^2f((k-2)x)-k^2f(kx)+2(k^2-1)f&((k-1)x)\|\\&\leq
\varphi((k-1)x,x) \hspace {0.75cm}(2.22)
\end{align*}
for all $x\in G.$ Replacing $x$ and $y$ by $(2k+1)x$ and $x$ in
$(2.15),$ respectively, we obtain
\begin{align*}
\|f((3k+1)x)-f(-(k+1)x)-k^2f(2(k+1)x)-k^2f(2kx)+&2(k^2-1)f((2k+1)x)\|\\&\leq
\varphi((2k+1)x,x) \hspace {.6cm}(2.23)
\end{align*}
for all $x\in G.$ Replacing $x$ and $y$ by $(2k-1)x$ and $x$ in
$(2.15),$ respectively, we have
\begin{align*}
\|f((3k-1)x)-f(-(k-1)x)-k^2f(2(k-1)x)-k^2f(2kx)+&2(k^2-1)f((2k-1)x)\|\\&\leq
\varphi((2k-1)x,x) \hspace {.6cm}(2.24)
\end{align*}
for all $x\in G.$ It follows from $(2.16),$ $(2.18),$ $(2.19),$
$(2.21)$ and $(2.22)$ that
\begin{align*}
\|f(3x)-4f(2x)+5f(x)\|\leq
&\frac{1}{|k^2(k^2-1)|}~\max\{\max\{~|2(k^2-1)|~\varphi(x,x),~|k^2|~\varphi(2x,x)\}\\&,\max\{\varphi(x,2x)
,\max\{\varphi((k+1)x,x),\varphi((k-1)x,x)\}~\}\} \hspace
{.5cm}(2.25)
\end{align*}
for all $x\in G.$ And, from $(2.16),$ $(2.17),$ $(2.19),$
$(2.20),$ $(2.23)$ and $(2.24),$ we conclude that
\begin{align*}
\|f&(4x)-2f(3x)-2f(2x)+6f(x)\|\leq
\frac{1}{|k^2(k^2-1)|}~\max\{\max\{\varphi(x,x),|k^2|\varphi(2x,2x)\}\\&,\max\{\max\{|2(k^2-1)|\varphi(x,2x)
,\varphi(x,3x)\},\max\{\varphi((2k+1)x,x),\varphi((2k-1)x,x)\}\}~\}
\hspace {.25cm}(2.26)
\end{align*}
for all $x\in G.$ Finally, by using $(2.25)$ and $(2.26),$ we
obtain that
\begin{align*}
\|f(4x)-10f(2x)&+16f(x)\|\leq
\frac{1}{|k^2(k^2-1)|}\max\{~\max\{\max\{|2(k^2-1)|\varphi(x,x),|k^2|\varphi(2x,x)\}\\&,\max\{\varphi(x,2x)
,\max\{\varphi((k+1)x,x),\varphi((k-1)x,x)\}~\}\},\max\{\max\{\varphi(x,x),|k^2|\varphi(2x,2x)\}\\&,\max\{\max\{|2(k^2-1)|\varphi(x,2x)
,\varphi(x,3x)\},\max\{\varphi((2k+1)x,x),\varphi((2k-1)x,x)\}\}~\}\}\\
&\hspace {9.75cm}(2.27)
\end{align*}
for all $x\in G.$ Let $g: G \to X$ be a function defined by
$g(x):=f(2x)-8f(x)$ for all $x \in G.$ From $(2.27),$ we conclude
that
\begin{align*}
\|\frac{g(2x)}{2}&-g(x)\|\leq
\frac{1}{|2.k^2(k^2-1)|}\max\{~\max\{\max\{|2(k^2-1)|\varphi(x,x),|k^2|\varphi(2x,x)\}\\&,\max\{\varphi(x,2x)
,\max\{\varphi((k+1)x,x),\varphi((k-1)x,x)\}~\}\},\max\{\max\{\varphi(x,x),|k^2|\varphi(2x,2x)\}\\&,\max\{\max\{|2(k^2-1)|\varphi(x,2x)
,\varphi(x,3x)\},\max\{\varphi((2k+1)x,x),\varphi((2k-1)x,x)\}\}~\}\}\\
&\hspace {11.25cm}(2.28)
\end{align*}
for all $x \in G.$ Replacing $x$ by $2^{n-1}x$ in $(2.28),$ we get
\begin{align*}
&\|\frac{g(2^nx)}{2^n}-\frac{g(2^{n-1}x)}{2^{n-1}}\|\leq
\frac{1}{|2^n.k^2(k^2-1)|}\\
&\max\{~\max\{\max\{|2(k^2-1)|\varphi(2^{n-1}x,2^{n-1}x),|k^2|\varphi(2^nx,2^{n-1}x)\}\\&,\max\{\varphi(2^{n-1}x,2^nx)
,\max\{\varphi(2^{n-1}(k+1)x,2^{n-1}x),\varphi(2^{n-1}(k-1)x,2^{n-1}x)\}~\}\}\\
&,\max\{\max\{\varphi(2^{n-1}x,2^{n-1}x),|k^2|\varphi(2^nx,2^nx)\}\\&,\max\{\max\{|2(k^2-1)|\varphi(2^{n-1}x,2^nx)
,\varphi(2^{n-1}x,3.2^{n-1}x)\}\\
&,\max\{\varphi(2^{n-1}(2k+1)x,2^{n-1}x),\varphi(2^{n-1}(2k-1)x,2^{n-1}x)\}\}~\}\}\\
&\hspace {11.25cm}(2.29)
\end{align*}
for all $x \in G.$ It follows from $(2.11)$ and $(2.29)$ that the
sequence $\{\frac{g(2^nx)}{2^{n}}\}$ is Cauchy. Since $X$ is
complete, we conclude that $\{\frac{g(2^nx)}{2^{n}}\}$ is
convergent. Set
$A(x):=\lim_{n \to \infty} \frac{g(2^nx)}{2^{n}}.$\\
Using induction one can show that
\begin{align*}
\|\frac{g(2^nx)}{2^{n}}-g(x)\| &\leq  \frac{1}{|2.k^2(k^2-1)|} \max \{\frac{1}{|2^i|}\\
&\max\{~\max\{\max\{|2(k^2-1)|\varphi(2^{i-1}x,2^{i-1}x),|k^2|\varphi(2^ix,2^{i-1}x)\}\\&,\max\{\varphi(2^{i-1}x,2^ix)
,\max\{\varphi(2^{i-1}(k+1)x,2^{i-1}x),\varphi(2^{i-1}(k-1)x,2^{i-1}x)\}~\}\}\\
&,\max\{\max\{\varphi(2^{i-1}x,2^{i-1}x),|k^2|\varphi(2^ix,2^ix)\}\\&,\max\{\max\{|2(k^2-1)|\varphi(2^{i-1}x,2^ix)
,\varphi(2^{i-1}x,3.2^{i-1}x)\}\\
&,\max\{\varphi(2^{i-1}(2k+1)x,2^{i-1}x),\varphi(2^{i-1}(2k-1)x,2^{i-1}x)\}\}~\}\}:~
0 \leq i < n \} \\
&\hspace{9cm}(2.30)
\end{align*}
 for all $n \in \Bbb N$ and all $x \in
G.$ By taking $n$ to approach infinity in $(2.30)$ and using
$(2.12)$ one obtains $(2.14).$ By $(2.10)$ and $(2.13),$ we get
\begin{align*}
\|DA(x,y)\|=\lim_{n \to \infty}\frac{1}{|2^{n}|}\|g(2^nx,2^ny)\|&=\lim_{n \to \infty}\frac{1}{|2^n|}\|D(f(2^{n+1}x,2^{n+1}y)-8f(2^nx,2^ny))\|\\
&\leq \lim_{n \to \infty}
\frac{1}{|2^n|}\max\{\|D(f(2^{n+1}x,2^{n+1}y)\|,\|8f(2^nx,2^ny))\|\}\\
&\leq\lim_{n \to \infty}\frac{1}{|2^n|}\max\{\varphi(2^{n+1}x,2^{n+1}y),|8|\varphi(2^nx,2^ny)\}\\
&=0
\end{align*}
for all $x,y \in G.$ Therefore the function $A:G \to X$ satisfies
$(1.3).$ If $A^{'}$ is another additive function satisfying
$(2.14),$ then
\begin{align*}
\|A(x)-A^{'}(x)\|&=\lim_{i \to
\infty}|2|^{-i}\|A(2^ix)-A^{'}(2^ix)\|\\
& \leq \lim_{i \to \infty} |2|^{-i} \max \{~
\|A(2^ix)-g(2^ix)\|,\|g(2^ix)-A^{'}(2^ix)\|~ \}\\
&\leq \frac{1}{|2.k^2(k^2-1)|} \lim_{i \to \infty} \lim_{n \to
\infty} \max \{\frac{1}{|2^j|}\\
&\max\{~\max\{\max\{|2(k^2-1)|\varphi(2^{j-1}x,2^{j-1}x),|k^2|\varphi(2^jx,2^{j-1}x)\}\\&,\max\{\varphi(2^{j-1}x,2^jx)
,\max\{\varphi(2^{j-1}(k+1)x,2^{j-1}x),\varphi(2^{j-1}(k-1)x,2^{j-1}x)\}~\}\}\\
&,\max\{\max\{\varphi(2^{j-1}x,2^{j-1}x),|k^2|\varphi(2^jx,2^jx)\}\\&,\max\{\max\{|2(k^2-1)|\varphi(2^{j-1}x,2^jx)
,\varphi(2^{j-1}x,3.2^{j-1}x)\}\\
&,\max\{\varphi(2^{j-1}(2k+1)x,2^{j-1}x),\varphi(2^{j-1}(2k-1)x,2^{j-1}x)\}\}~\}\}:~
i \leq j < n+i \}\\
&=0
\end{align*}
for all $x \in G.$ Therefore $A=A^{'}.$ This completes the proof
of the uniqueness of $A.$
\end{proof}

\begin{thm}
Let $\varphi:G \times G \to [0,\infty)$ be a function such that
$$\lim_{n \to \infty} \frac{1}{|8^n|}\max\{\varphi(2^{n+1}x,2^{n+1}y),|2|\varphi(2^nx,2^ny)\}=0\eqno(2.31)$$
\begin{align*}
\lim_{n \to \infty} \frac{1}{|8^n.k^2(k^2-1)|}
&\max\{~\max\{\max\{|2(k^2-1)|\varphi(2^{n-1}x,2^{n-1}x),|k^2|\varphi(2^nx,2^{n-1}x)\}\\&,\max\{\varphi(2^{n-1}x,2^nx)
,\max\{\varphi(2^{n-1}(k+1)x,2^{n-1}x),\varphi(2^{n-1}(k-1)x,2^{n-1}x)\}~\}\}\\
&,\max\{\max\{\varphi(2^{n-1}x,2^{n-1}x),|k^2|\varphi(2^nx,2^nx)\}\\&,\max\{\max\{|2(k^2-1)|\varphi(2^{n-1}x,2^nx)
,\varphi(2^{n-1}x,3.2^{n-1}x)\}\\
&,\max\{\varphi(2^{n-1}(2k+1)x,2^{n-1}x),\varphi(2^{n-1}(2k-1)x,2^{n-1}x)\}\}~\}\}=0\\
&\hspace{9cm}(2.32)
\end{align*}
for all $x,y \in G$ and let for each $x \in G$ the limit
\begin{align*}
\max
&\{\frac{1}{|8^i|}\max\{~\max\{\max\{|2(k^2-1)|\varphi(2^{i-1}x,2^{i-1}x),|k^2|\varphi(2^ix,2^{i-1}x)\}\\&,\max\{\varphi(2^{i-1}x,2^ix)
,\max\{\varphi(2^{i-1}(k+1)x,2^{i-1}x),\varphi(2^{i-1}(k-1)x,2^{i-1}x)\}~\}\}\\
&,\max\{\max\{\varphi(2^{i-1}x,2^{i-1}x),|k^2|\varphi(2^ix,2^ix)\}\\&,\max\{\max\{|2(k^2-1)|\varphi(2^{i-1}x,2^ix)
,\varphi(2^{i-1}x,3.2^{i-1}x)\}\\
&,\max\{\varphi(2^{i-1}(2k+1)x,2^{i-1}x),\varphi(2^{i-1}(2k-1)x,2^{i-1}x)\}\}~\}\}:~
0 \leq i < n \},\hspace{1cm}(2.33)
\end{align*}
denoted by $\tilde{\varphi}_C(x),$ exist. Suppose that $f:G \to X$
is an odd function satisfying
$$\|Df(x,y)\| \leq \varphi(x,y) \eqno(2.34)$$
for all $x,y \in G.$ Then there exist a cubic function $C:G \to X$
such that
$$\|f(2x)-2f(x)-C(x)\| \leq \frac{1}{|8.k^2(k^2-1)|} \tilde{\varphi}_C(x) \eqno(2.35)$$
for all $x \in G.$ Moreover, if
\begin{align*}
\lim_{i \to \infty} \lim_{n \to \infty} \max
&\{\frac{1}{|8^j|}\max\{~\max\{\max\{|2(k^2-1)|\varphi(2^{j-1}x,2^{j-1}x),|k^2|\varphi(2^jx,2^{j-1}x)\}\\&,\max\{\varphi(2^{j-1}x,2^jx)
,\max\{\varphi(2^{j-1}(k+1)x,2^{j-1}x),\varphi(2^{j-1}(k-1)x,2^{j-1}x)\}~\}\}\\
&,\max\{\max\{\varphi(2^{j-1}x,2^{j-1}x),|k^2|\varphi(2^jx,2^jx)\}\\&,\max\{\max\{|2(k^2-1)|\varphi(2^{j-1}x,2^jx)
,\varphi(2^{j-1}x,3.2^{j-1}x)\}\\
&,\max\{\varphi(2^{j-1}(2k+1)x,2^{j-1}x),\varphi(2^{j-1}(2k-1)x,2^{j-1}x)\}\}~\}\}:~
i \leq j < n+i \}\\
&=0,
\end{align*}
then $C$ is the unique cubic function satisfying $(2.35).$
\end{thm}
\begin{proof}
Similar to the proof of Theorem 2.2, we have
\begin{align*}
\|f(4x)-10f(2x)&+16f(x)\|\leq
\frac{1}{|k^2(k^2-1)|}\max\{~\max\{\max\{|2(k^2-1)|\varphi(x,x),|k^2|\varphi(2x,x)\}\\&,\max\{\varphi(x,2x)
,\max\{\varphi((k+1)x,x),\varphi((k-1)x,x)\}~\}\},\max\{\max\{\varphi(x,x),|k^2|\varphi(2x,2x)\}\\&,\max\{\max\{|2(k^2-1)|\varphi(x,2x)
,\varphi(x,3x)\},\max\{\varphi((2k+1)x,x),\varphi((2k-1)x,x)\}\}~\}\}
\end{align*}
for all $x\in G.$ Let $h: G \to X$ be a function defined by
$h(x):=f(2x)-2f(x)$ for all $x \in G$ then we have
\begin{align*}
\|\frac{h(2x)}{8}&-h(x)\|\leq
\frac{1}{|8.k^2(k^2-1)|}\max\{~\max\{\max\{|2(k^2-1)|\varphi(x,x),|k^2|\varphi(2x,x)\}\\&,\max\{\varphi(x,2x)
,\max\{\varphi((k+1)x,x),\varphi((k-1)x,x)\}~\}\},\max\{\max\{\varphi(x,x),|k^2|\varphi(2x,2x)\}\\&,\max\{\max\{|2(k^2-1)|\varphi(x,2x)
,\varphi(x,3x)\},\max\{\varphi((2k+1)x,x),\varphi((2k-1)x,x)\}\}~\}\}\\
&\hspace {11.25cm}(2.36)
\end{align*}
for all $x \in G.$ Replacing $x$ by $2^{n-1}x$ in $(2.36),$ we get
\begin{align*}
&\|\frac{h(2^nx)}{8^n}-\frac{h(2^{n-1}x)}{8^{n-1}}\|\leq
\frac{1}{|8^n.k^2(k^2-1)|}\\
&\max\{~\max\{\max\{|2(k^2-1)|\varphi(2^{n-1}x,2^{n-1}x),|k^2|\varphi(2^nx,2^{n-1}x)\}\\&,\max\{\varphi(2^{n-1}x,2^nx)
,\max\{\varphi(2^{n-1}(k+1)x,2^{n-1}x),\varphi(2^{n-1}(k-1)x,2^{n-1}x)\}~\}\}\\
&,\max\{\max\{\varphi(2^{n-1}x,2^{n-1}x),|k^2|\varphi(2^nx,2^nx)\}\\&,\max\{\max\{|2(k^2-1)|\varphi(2^{n-1}x,2^nx)
,\varphi(2^{n-1}x,3.2^{n-1}x)\}\\
&,\max\{\varphi(2^{n-1}(2k+1)x,2^{n-1}x),\varphi(2^{n-1}(2k-1)x,2^{n-1}x)\}\}~\}\}\hspace
{3.75cm}(2.37)
\end{align*}
for all $x \in G.$ It follows from $(2.32)$ and $(2.37)$ that the
sequence $\{\frac{h(2^nx)}{8^{n}}\}$ is Cauchy. Since $X$ is
complete, we conclude that $\{\frac{h(2^nx)}{8^{n}}\}$ is
convergent. Set
$C(x):=\lim_{n \to \infty} \frac{h(2^nx)}{8^{n}}.$\\
Using induction one can show that
\begin{align*}
\|\frac{h(2^nx)}{8^{n}}-h(x)\| &\leq  \frac{1}{|8.k^2(k^2-1)|} \max \{\frac{1}{|8^i|}\\
&\max\{~\max\{\max\{|2(k^2-1)|\varphi(2^{i-1}x,2^{i-1}x),|k^2|\varphi(2^ix,2^{i-1}x)\}\\&,\max\{\varphi(2^{i-1}x,2^ix)
,\max\{\varphi(2^{i-1}(k+1)x,2^{i-1}x),\varphi(2^{i-1}(k-1)x,2^{i-1}x)\}~\}\}\\
&,\max\{\max\{\varphi(2^{i-1}x,2^{i-1}x),|k^2|\varphi(2^ix,2^ix)\}\\&,\max\{\max\{|2(k^2-1)|\varphi(2^{i-1}x,2^ix)
,\varphi(2^{i-1}x,3.2^{i-1}x)\}\\
&,\max\{\varphi(2^{i-1}(2k+1)x,2^{i-1}x),\varphi(2^{i-1}(2k-1)x,2^{i-1}x)\}\}~\}\}:~
0 \leq i < n \} \\
&\hspace{9cm}(2.38)
\end{align*}
 for all $n \in \Bbb N$ and all $x \in
G.$ By taking $n$ to approach infinity in $(2.38)$ and using
$(2.33)$ one obtains $(2.35).$ By $(2.31)$ and $(2.34),$ we get
\begin{align*}
\|DC(x,y)\|=\lim_{n \to \infty}\frac{1}{|8^{n}|}\|h(2^nx,2^ny)\|&=\lim_{n \to \infty}\frac{1}{|8^n|}\|D(f(2^{n+1}x,2^{n+1}y)-2f(2^nx,2^ny))\|\\
&\leq \lim_{n \to \infty}
\frac{1}{|8^n|}\max\{\|D(f(2^{n+1}x,2^{n+1}y)\|,\|2f(2^nx,2^ny))\|\}\\
&\leq\lim_{n \to \infty}\frac{1}{|8^n|}\max\{\varphi(2^{n+1}x,2^{n+1}y),|2|\varphi(2^nx,2^ny)\}\\
&=0
\end{align*}
for all $x,y \in G.$ Therefore the function $C:G \to X$ satisfies
$(1.3).$ If $C^{'}$ is another cubic function satisfying $(2.35),$
then
\begin{align*}
\|C(x)-C^{'}(x)\|&=\lim_{i \to
\infty}|8|^{-i}\|C(2^ix)-C^{'}(2^ix)\|\\
& \leq \lim_{i \to \infty} |8|^{-i} \max \{~
\|C(2^ix)-h(2^ix)\|,\|h(2^ix)-C^{'}(2^ix)\|~ \}\\
&\leq \frac{1}{|8.k^2(k^2-1)|} \lim_{i \to \infty} \lim_{n \to
\infty} \max \{\frac{1}{|8^j|}\\
&\max\{~\max\{\max\{|2(k^2-1)|\varphi(2^{j-1}x,2^{j-1}x),|k^2|\varphi(2^jx,2^{j-1}x)\}\\&,\max\{\varphi(2^{j-1}x,2^jx)
,\max\{\varphi(2^{j-1}(k+1)x,2^{j-1}x),\varphi(2^{j-1}(k-1)x,2^{j-1}x)\}~\}\}\\
&,\max\{\max\{\varphi(2^{j-1}x,2^{j-1}x),|k^2|\varphi(2^jx,2^jx)\}\\&,\max\{\max\{|2(k^2-1)|\varphi(2^{j-1}x,2^jx)
,\varphi(2^{j-1}x,3.2^{j-1}x)\}\\
&,\max\{\varphi(2^{j-1}(2k+1)x,2^{j-1}x),\varphi(2^{j-1}(2k-1)x,2^{j-1}x)\}\}~\}\}:~
i \leq j < n+i \}\\
&=0
\end{align*}
for all $x \in G.$ Therefore $C=C^{'}.$ This completes the proof
of the uniqueness of $C.$
\end{proof}

\begin{thm}
Let $\varphi:G \times G \to [0,\infty)$ be a function such that
\begin{align*}
&\lim_{n \to \infty}
\frac{1}{|2^n|}\max\{\varphi(2^{n+1}x,2^{n+1}y),|8|\varphi(2^nx,2^ny)\}\\
&=\lim_{n \to \infty}
\frac{1}{|8^n|}\max\{\varphi(2^{n+1}x,2^{n+1}y),|2|\varphi(2^nx,2^ny)\}=0
\hspace{3.5cm}(2.39)
\end{align*}
for all $x,y \in G$ and let for each $x \in G$ the limit
\begin{align*}
\max
&\{\frac{1}{|2^i|}\max\{~\max\{\max\{|2(k^2-1)|\varphi(2^{i-1}x,2^{i-1}x),|k^2|\varphi(2^ix,2^{i-1}x)\}\\&,\max\{\varphi(2^{i-1}x,2^ix)
,\max\{\varphi(2^{i-1}(k+1)x,2^{i-1}x),\varphi(2^{i-1}(k-1)x,2^{i-1}x)\}~\}\}\\
&,\max\{\max\{\varphi(2^{i-1}x,2^{i-1}x),|k^2|\varphi(2^ix,2^ix)\}\\&,\max\{\max\{|2(k^2-1)|\varphi(2^{i-1}x,2^ix)
,\varphi(2^{i-1}x,3.2^{i-1}x)\}\\
&,\max\{\varphi(2^{i-1}(2k+1)x,2^{i-1}x),\varphi(2^{i-1}(2k-1)x,2^{i-1}x)\}\}~\}\}:~
0 \leq i < n \},
\end{align*}
denoted by $\tilde{\varphi}_A(x),$ and
\begin{align*}
\max
&\{\frac{1}{|8^i|}\max\{~\max\{\max\{|2(k^2-1)|\varphi(2^{i-1}x,2^{i-1}x),|k^2|\varphi(2^ix,2^{i-1}x)\}\\&,\max\{\varphi(2^{i-1}x,2^ix)
,\max\{\varphi(2^{i-1}(k+1)x,2^{i-1}x),\varphi(2^{i-1}(k-1)x,2^{i-1}x)\}~\}\}\\
&,\max\{\max\{\varphi(2^{i-1}x,2^{i-1}x),|k^2|\varphi(2^ix,2^ix)\}\\&,\max\{\max\{|2(k^2-1)|\varphi(2^{i-1}x,2^ix)
,\varphi(2^{i-1}x,3.2^{i-1}x)\}\\
&,\max\{\varphi(2^{i-1}(2k+1)x,2^{i-1}x),\varphi(2^{i-1}(2k-1)x,2^{i-1}x)\}\}~\}\}:~
0 \leq i < n \}
\end{align*}
denoted by $\tilde{\varphi}_C(x),$ exists. Suppose that $f:G \to
X$ is an odd function satisfying
$$\|Df(x,y)\| \leq \varphi(x,y)\eqno(2.40)$$
for all $x,y \in G.$ Then there exists an additive function $A:G
\to X$ and a cubic function $C:G \to X$ such that
$$\|f(x)-A(x)-C(x)\| \leq \frac{1}{|k^2(k^2-1)|}\max
\{\frac{1}{|2|} \tilde{\varphi}_A(x),\frac{1}{|8|}
\tilde{\varphi}_C(x) \}\eqno(2.41)$$ for all $x \in G.$ Moreover,
if
\begin{align*}
\lim_{i \to \infty} \lim_{n \to \infty} \max
&\{\frac{1}{|2^j|}\max\{~\max\{\max\{|2(k^2-1)|\varphi(2^{j-1}x,2^{j-1}x),|k^2|\varphi(2^jx,2^{j-1}x)\}\\&,\max\{\varphi(2^{j-1}x,2^jx)
,\max\{\varphi(2^{j-1}(k+1)x,2^{j-1}x),\varphi(2^{j-1}(k-1)x,2^{j-1}x)\}~\}\}\\
&,\max\{\max\{\varphi(2^{j-1}x,2^{j-1}x),|k^2|\varphi(2^jx,2^jx)\}\\&,\max\{\max\{|2(k^2-1)|\varphi(2^{j-1}x,2^jx)
,\varphi(2^{j-1}x,3.2^{j-1}x)\}\\
&,\max\{\varphi(2^{j-1}(2k+1)x,2^{j-1}x),\varphi(2^{j-1}(2k-1)x,2^{j-1}x)\}\}~\}\}:~
i \leq j < n+i \}\\
&=0,
\end{align*}
\begin{align*}
\lim_{i \to \infty} \lim_{n \to \infty} \max
&\{\frac{1}{|8^j|}\max\{~\max\{\max\{|2(k^2-1)|\varphi(2^{j-1}x,2^{j-1}x),|k^2|\varphi(2^jx,2^{j-1}x)\}\\&,\max\{\varphi(2^{j-1}x,2^jx)
,\max\{\varphi(2^{j-1}(k+1)x,2^{j-1}x),\varphi(2^{j-1}(k-1)x,2^{j-1}x)\}~\}\}\\
&,\max\{\max\{\varphi(2^{j-1}x,2^{j-1}x),|k^2|\varphi(2^jx,2^jx)\}\\&,\max\{\max\{|2(k^2-1)|\varphi(2^{j-1}x,2^jx)
,\varphi(2^{j-1}x,3.2^{j-1}x)\}\\
&,\max\{\varphi(2^{j-1}(2k+1)x,2^{j-1}x),\varphi(2^{j-1}(2k-1)x,2^{j-1}x)\}\}~\}\}:~
i \leq j < n+i \}\\
&=0,
\end{align*}
then $A$ is the unique additive function and $C$ is the unique
cubic function satisfying $(2.41).$
\end{thm}

\begin{proof}
By Theorems 2.2 and 2.3, there exists an additive function $A_1:G
\to X$ and a cubic function $C_1:G \to X$ such that
$$\|f(2x)-8f(x)-A(x)\| \leq \frac{1}{|2.k^2(k^2-1)|} \tilde{\varphi}_A(x)$$
$$\|f(2x)-2f(x)-C(x)\| \leq \frac{1}{|8.k^2(k^2-1)|} \tilde{\varphi}_C(x)$$
for all $x \in G.$ So we obtain $(2.41)$ by letting
$A(x)=\frac{-1}{6} A_1(x)$ and $C(x)=\frac{1}{6}C_1(x)$ for all $x
\in G.$ Now it is obvious that $(2.41)$ holds true for all $x \in
G,$ and the proof of theorem is complete.
\end{proof}

\begin{thm}
Let $\varphi:G \times G \to [0,\infty)$ be a function such that
\begin{align*}
&\lim_{n \to \infty}
\frac{1}{|2^n|}\max\{\varphi(2^{n+1}x,2^{n+1}y),|8|\varphi(2^nx,2^ny)\}\\
&=\lim_{n \to \infty}
\frac{1}{|8^n|}\max\{\varphi(2^{n+1}x,2^{n+1}y),|2|\varphi(2^nx,2^ny)\}\\
&=\lim_{n \to \infty} \frac{\varphi(2^nx,2^ny)}{|2|^{2n}}=0
\hspace{8.75cm}(2.42)
\end{align*}
for all $x,y \in G$ and let for each $x \in G$ the limit
\begin{align*}
\max
&\{\frac{1}{|2^i|}\max\{~\max\{\max\{|2(k^2-1)|\varphi(2^{i-1}x,2^{i-1}x),|k^2|\varphi(2^ix,2^{i-1}x)\}\\&,\max\{\varphi(2^{i-1}x,2^ix)
,\max\{\varphi(2^{i-1}(k+1)x,2^{i-1}x),\varphi(2^{i-1}(k-1)x,2^{i-1}x)\}~\}\}\\
&,\max\{\max\{\varphi(2^{i-1}x,2^{i-1}x),|k^2|\varphi(2^ix,2^ix)\}\\&,\max\{\max\{|2(k^2-1)|\varphi(2^{i-1}x,2^ix)
,\varphi(2^{i-1}x,3.2^{i-1}x)\}\\
&,\max\{\varphi(2^{i-1}(2k+1)x,2^{i-1}x),\varphi(2^{i-1}(2k-1)x,2^{i-1}x)\}\}~\}\}:~
0 \leq i < n \},
\end{align*}
denoted by $\tilde{\varphi}_A(x),$ and
$$\lim_{n \to \infty} \max \{\frac{1}{|k^{2j}|}\varphi(0,k^jx) :~ 0 \leq j < n \},$$
denoted by $\tilde{\varphi}_Q(x),$ and
\begin{align*}
\max
&\{\frac{1}{|8^i|}\max\{~\max\{\max\{|2(k^2-1)|\varphi(2^{i-1}x,2^{i-1}x),|k^2|\varphi(2^ix,2^{i-1}x)\}\\&,\max\{\varphi(2^{i-1}x,2^ix)
,\max\{\varphi(2^{i-1}(k+1)x,2^{i-1}x),\varphi(2^{i-1}(k-1)x,2^{i-1}x)\}~\}\}\\
&,\max\{\max\{\varphi(2^{i-1}x,2^{i-1}x),|k^2|\varphi(2^ix,2^ix)\}\\&,\max\{\max\{|2(k^2-1)|\varphi(2^{i-1}x,2^ix)
,\varphi(2^{i-1}x,3.2^{i-1}x)\}\\
&,\max\{\varphi(2^{i-1}(2k+1)x,2^{i-1}x),\varphi(2^{i-1}(2k-1)x,2^{i-1}x)\}\}~\}\}:~
0 \leq i < n \}
\end{align*}
denoted by $\tilde{\varphi}_C(x),$ exists. Suppose that $f:G \to
X$ is a function satisfying
$$\|Df(x,y)\| \leq \varphi(x,y)\eqno(2.43)$$
for all $x,y \in G.$ Then there exists an additive function $A:G
\to X$ and a quadratic function $Q:G \to X$ and a cubic function
$C:G \to X$ such that
\begin{align*}
\|f(x)-A(x)-Q(x)-C(x)\|& \leq
\max\{\frac{1}{|2k^2(k^2-1)|}\max\{\max \{\frac{1}{|2|}
\tilde{\varphi}_A(x),\frac{1}{|8|} \tilde{\varphi}_C(x) \}\\
&,\max \{\frac{1}{|2|} \tilde{\varphi}_A(-x),\frac{1}{|8|}
\tilde{\varphi}_C(-x) \}\},\frac{1}{|4.k^2|}
\max\{\tilde{\varphi}_Q(x),\tilde{\varphi}_Q(x)\}\}\\
&\hspace{7.5cm}(2.44)
\end{align*}
for all $x \in G.$ Moreover, if
\begin{align*}
\lim_{i \to \infty} \lim_{n \to \infty} \max
&\{\frac{1}{|2^j|}\max\{~\max\{\max\{|2(k^2-1)|\varphi(2^{j-1}x,2^{j-1}x),|k^2|\varphi(2^jx,2^{j-1}x)\}\\&,\max\{\varphi(2^{j-1}x,2^jx)
,\max\{\varphi(2^{j-1}(k+1)x,2^{j-1}x),\varphi(2^{j-1}(k-1)x,2^{j-1}x)\}~\}\}\\
&,\max\{\max\{\varphi(2^{j-1}x,2^{j-1}x),|k^2|\varphi(2^jx,2^jx)\}\\&,\max\{\max\{|2(k^2-1)|\varphi(2^{j-1}x,2^jx)
,\varphi(2^{j-1}x,3.2^{j-1}x)\}\\
&,\max\{\varphi(2^{j-1}(2k+1)x,2^{j-1}x),\varphi(2^{j-1}(2k-1)x,2^{j-1}x)\}\}~\}\}:~
i \leq j < n+i \}\\
&=0,
\end{align*}
\begin{align*}
\lim_{i \to \infty} \lim_{n \to \infty}
 \max \{\max \{\frac{1}{|k^{2j}|}\varphi(0,k^jx)
\}:~i \leq j < n+i \}=0,
\end{align*}
\begin{align*}
\lim_{i \to \infty} \lim_{n \to \infty} \max
&\{\frac{1}{|8^j|}\max\{~\max\{\max\{|2(k^2-1)|\varphi(2^{j-1}x,2^{j-1}x),|k^2|\varphi(2^jx,2^{j-1}x)\}\\&,\max\{\varphi(2^{j-1}x,2^jx)
,\max\{\varphi(2^{j-1}(k+1)x,2^{j-1}x),\varphi(2^{j-1}(k-1)x,2^{j-1}x)\}~\}\}\\
&,\max\{\max\{\varphi(2^{j-1}x,2^{j-1}x),|k^2|\varphi(2^jx,2^jx)\}\\&,\max\{\max\{|2(k^2-1)|\varphi(2^{j-1}x,2^jx)
,\varphi(2^{j-1}x,3.2^{j-1}x)\}\\
&,\max\{\varphi(2^{j-1}(2k+1)x,2^{j-1}x),\varphi(2^{j-1}(2k-1)x,2^{j-1}x)\}\}~\}\}:~
i \leq j < n+i \}\\
&=0,
\end{align*}
then $A$ is the unique additive function and $Q$ is the unique
quadratic function and $C$ is the unique cubic function.
\end{thm}

\begin{proof}
Let $f_o(x)=\frac{1}{2}[f(x)-f(-x)]$ for all $x\in G.$ Then
$f_o(0)=0,$ $f_o(-x)=-f_o(x),$ and $$\|Df_o(x,y)\|\leq
\frac{1}{|2|}\max\{\varphi(x,y),\varphi(-x,-y)\}$$ \hspace {0.3cm}
for all $x,y\in G.$ From Theorem 2.3, it follows that there exists
a unique additive function $A:G \to X$ and a unique cubic function
$C:G \to X$ satisfying
$$\|f_o(x)-A(x)-C(x)\| \leq \frac{1}{|2k^2(k^2-1)|}\max\{\max
\{\frac{1}{|2|} \tilde{\varphi}_A(x),\frac{1}{|8|}
\tilde{\varphi}_C(x) \},\max \{\frac{1}{|2|}
\tilde{\varphi}_A(-x),\frac{1}{|8|} \tilde{\varphi}_C(-x)
\}\}\eqno(2.45)$$
for all $x\in G.$\\
Let $f_e(x)=\frac{1}{2}[f(x)+f(-x)]$ for all $x\in G.$ Then
$f_e(0)=0,$ $f_e(-x)=f_e(x),$ and $$\|Df_e(x,y)\| \leq
\frac{1}{|2|}\max\{\varphi(x,y),\varphi(-x,-y)\}$$ \hspace {0.3cm}
for all $x,y\in G.$ From Theorem 2.1, it follows that there exist
a unique quadratic function $Q:G \to X$ satisfying
$$\|f_e(x)-Q(x)\| \leq \frac{1}{|4.k^2|} \max\{\tilde{\varphi}_Q(x),\tilde{\varphi}_Q(x)\} \eqno(2.46)$$
for all $x\in G.$ Hence $(2.44)$ follows from $(2.45)$ and
$(2.46).$
\end{proof}

{\small


}
\end{document}